\documentclass[11pt]{amsart}
\usepackage{amssymb,latexsym,mathtools}
\usepackage{hyperref}
\hypersetup{
    pdftoolbar=true,        
    pdfmenubar=false,        
    pdffitwindow=false,     
    pdfstartview={FitH},    % fits the width of the page to the window
    pdftitle={},
    pdfauthor={},
    pdfsubject={},
    pdfkeywords={},
    pdfnewwindow=true,      % links in new window
    colorlinks=true,       % false: boxed links; true: colored links
    linkcolor=red,          % color of internal links
    citecolor=blue,        % color of links to bibliography
    urlcolor=cyan,           % color of external links
    }
\newtheorem{theorem}{Theorem}[section]
\newtheorem{proposition}[theorem]{Proposition}

\newtheorem{lemma}[theorem]{Lemma}
\newtheorem{remark}[theorem]{Remark}

\newcommand{\des}{{\rm des}}
\newcommand{\Des}{{\rm Des}}
\newcommand{\SYT}{{\rm SYT}}
\newcommand{\SYB}{{\rm SYB}}
\newcommand{\QSym}{{\rm QSym}}  
\newcommand{\Par}{{\rm Par}}
\newcommand{\sDes}{{\rm sDes}}
\newcommand{\NN}{{\mathbb N}}
\newcommand{\ZZ}{{\mathbb Z}}
\newcommand{\CC}{{\mathbb C}}

\begin{document}
\title[The Eulerian distribution on the involutions of the hyperoctahedral group is unimodal]{The Eulerian distribution on the involutions of the hyperoctahedral group is unimodal}

\author{Vassilis-Dionyssis~P.~Moustakas}

\address{Department of Mathematics,
National and Kapodistrian University of Athens\\
Panepistimioupolis\\
15784 Athens, Greece}
\email{vasmous@math.uoa.gr}

\date{January 22, 2018}
\thanks{ \textit{Key words and phrases}. 
Involution, Descent, Unimodality, Hyperoctahedral group, Quasisymmetric function.}

\begin{abstract} 
The Eulerian distribution on the involutions of the symmetric group is unimodal, as shown by Guo and Zeng. In this paper we prove that the Eulerian distribution on the involutions of the hyperoctahedral group, when viewed as a colored permutation group, is unimodal in a similar way and we compute its generating function, using signed quasisymmetric functions.
\end{abstract}

\maketitle

\section{Introduction and results}
\label{sec:intro}

An involution is a permutation $w \in \mathfrak{S}_n$ such that $w^{-1} = w$. When written in cycle notation, such a permutation consists only of one-cycles and two-cycles. Let $I_{n}$ be the set of all involutions in the symmetric group $\mathfrak{S}_n$ and let
\[ I_{n}(x) = \sum_{w \in I_{n}} x^{\des(w)},\]
where $\des(w)$ is the number of descents (see Section \ref{sec:back} for missing definitions) of $w \in \mathfrak{S}_n$. For the first few values of $n$ we have:

\[I_{n}(x) =\begin{cases}
 1, & \text{if }  n=1\\
 1 +x, & \text{if }  n=2\\
 1 + 2x + x^{2}, & \text{if }  n=3 \\
 1 + 4x + 4x^{2} + x^{3}, & \text{if }  n=4 \\
 1 + 6x + 12x^{2} + 6x^{3} + x^{4}, & \text{if }  n=5 \\
  1 + 9x + 28x^{2} + 28x^{3} + 9x^{4} + x^{5}, &\text{if } n=6.
\end{cases}\]

\medskip
\noindent
Strehl \cite{Str81} proved that $I_{n}(x)$ is symmetric with center of symmetry at $(n-1)/2$ (conjectured by D. Dumont). Dukes \cite{Duk07} proved partially that $I_{n}(x)$ is unimodal, which was later fully established by Guo and Zeng \cite{GZ06}.

We will be concerned with a natural analogue of $I_{n}(x)$ for the hyperoctahedral group $B_{n}$ of signed permutations. It is defined by the formula
\[ I_{n}^{B}(x) = \sum_{w \in I_{n}^{B}} x^{\des_{B}(w)}\]
where $\des_{B}(w)$ is the number of descents of $w \in B_{n}$, when viewed as a colored permutation (see Section \ref{subsec:hyper} for the definition) and $I_{n}^{B}$ is the set of involutions in $B_{n}$. For the first few values of $n$ we have:

\[I_{n}^{B}(x) =\begin{cases}
 1 +x, & \text{if } n=1 \\
 1 + 4x + x^{2}, & \text{if } n=2 \\
 1 + 9x + 9x^{2} + x^{3}, & \text{if } n=3 \\
 1 + 17x + 40x^{2} + 17x^{3} + x^{4}, & \text{if } n=4 \\
 1 + 28x + 127x^{2} + 127x^{3} + 28x^{4} + x^{5}, & \text{if } n=5 \\
 1 + 43x + 331x^{2} + 632x^{3} + 331x^{4}+ 43x^{5} + x^{6}, & \text{if } n=6.
\end{cases}\]

\medskip
\noindent
 D\'esarm\'enien and Foata \cite{DF85} and later, using different methods, Gessel and Reutenauer \cite{GR93} computed the generating function 
 
\begin{equation} \label{eq:genfuninvo}
\sum_{n \geq 0} \, \frac{I_{n}(x)}{(1-x)^{n+1}} \, t^{n} = \sum_{m \geq 0} \, \frac{x^{m}}{(1-t)^{m+1}(1-t^{2})^{\frac{m(m+1)}{2}}}  
\end{equation}

\medskip
\noindent
of $I_{n}(x)$. Our first result computes the generating function of $I_{n}^{B}(x)$, using methods similar to those of \cite{GR93}.

\begin{theorem}\label{thm:genfuninvoB}
\begin{equation}\label{eq:genfuninvoB}
\sum_{n\geq 0} \, \frac{I_{n}^{B}(x)}{(1-x)^{n+1}} \, t^{n} = \sum_{m \geq 0} \, \frac{x^{m}}{(1-t)^{2m+1}(1-t^{2})^{m^{2}}}.
\end{equation}
\end{theorem}

\medskip
One key ingredient of the proof of Theorem \ref{thm:genfuninvoB} is a $B_{n}$-analogue of the well known expansion of the Schur symmetric functions in terms of fundamental quasisymmetric functions, obtained by Adin et al. in \cite{AAER17} (discussed in Subsection \ref{subsec:hyper}).

Using the generating function obtained in Theorem \ref{thm:genfuninvoB} we derive a linear recurrence formula for the sequence of coefficients of $I_{n}^{B}(x)$, which we use to prove the main result of this paper, namely the unimodality of $I_{n}^{B}(x)$. It is noted that a consequence of this formula is the palindromicity of $I_{n}^{B}(x)$, which can also be proved combinatorially using the natural $B_{n}$-analogue of the Robinson-Schensted correspondence (discussed in Subsection \ref{subsec:hyper}).

The structure of the paper is as follows. Section 2 fixes notation and reviews background material. Section 3 proves Theorem \ref{thm:genfuninvoB} and obtains a linear recurrence formula for the coefficients of $I_{n}^{B}(x)$. Section 4 proves the unimodality of $I_{n}(x)$. Section 5 discusses some remarks and open problems.
\section{Background and notation}
\label{sec:back}

This section fixes notation and briefly reviews background material regarding the combinatorics of (signed) permutations and Young (bi)tableaux, symmetric and unimodal polynomials and the theory of symmetric and quasisymmetric functions which will be needed in the sequel. More information on these topics can be found in \cite{AAER17}, \cite{StaEC1} and \cite[Chapter~7]{StaEC2}.

For positive integer $n$ we set $[n] \coloneqq \{1,2,\dots, n\}$ and $\Omega_{n} \coloneqq \{1,-1,2,-2,\dots,n,-n\}$. We denote by $|S|$ the cardinality of a finite set $S$.

\subsection{Permutations, tableaux and unimodal polynomials}
\label{subsec:perm}

We will denote by $\mathfrak{S}_{n}$ the symmetric group of all permutations of the set $[n]$, i.e. bijective maps $w : [n] \rightarrow [n]$. If $w = w_{1}w_{2}\cdots w_{n} \in \mathfrak{S}_{n}$ and $1 \leq i \leq n-1$, then $i$ is a descent of $w,$ if $w_{i} > w_{i+1}$ (otherwise, it is an ascent). Let $\Des(w)$ be the set of descents of $w \in \mathfrak{S}_{n}$ and set $\des(w) = |\Des(w)|$. A statistic on $\mathfrak{S}_{n}$ is called Eulerian, if it is equidistributed with the descent number statistic. The $n$th Eulerian polynomial \cite[Section~1.4]{StaEC1} is defined by the formula
\[ A_{n}(x) = \sum_{w \in \mathfrak{S}_{n}} x^{\des(w)}\]
for every positive integer $n$. 

Let $p(x) = \sum_{k=0}^{d}a_{k}x^{k}$ be a polynomial with real coefficients. We recall that $p(x)$ is unimodal (and has unimodal coefficients) if there exists an index $0 \leq j \leq d$ such that $a_{0} \leq a_{1} \leq \cdots \leq a_{j} \geq a_{j+1} \geq \cdots \geq a_{d}$. The polynomial $p(x)$ is said to be log-concave (and has log-concave coefficients) if $a_{i}^{2} \geq a_{i-1}a_{i+1}$ for $1 \leq i \leq d-1$. We will say that $p(x)$ is symmetric (and that it has symmetric coefficients) if there exists an integer $n \geq d$ such that $a_{i} = a_{n-i}$ for $0 \leq i \leq n$, where $a_{k}=0$ for $n>d$. The center of symmetry of $p(x)$ is then defined to be $n/2$ (provided $p(x)$ in nonzero). We say that $p(x)$ is $\gamma$-positive if
\[ p(x) = \sum_{i=0}^{\lfloor n/2 \rfloor}\gamma_{i}x^{i} (1+x)^{n-2i}\]
for some $n \in \NN$ and nonnegative reals $\gamma_{0}, \gamma_{1}, \dots, \gamma_{\lfloor n/2 \rfloor}$. Every $\gamma$-positive polynomial is symmetric and  unimodal, as a sum of  symmetric and unimodal polynomials with a common center of symmetry. The $n$th Eulerian polynomial is an example of a $\gamma$-positive polynomial (hence unimodal) as proved by Foata and Sch{\"u}tzenberger \cite{FS70} (see \cite[Theorem~2.1]{Athtba} for several combinatorial interpretations of the corresponding $\gamma$-coefficients). Unimodal polynomials arise often in combinatorics, geometry, and algebra, see for example \cite{Bra15} and \cite{Sta89}. The $\gamma$-positivity of $p(x)$ often provides a more elementary proof of the unimodality of $p(x)$. For a comprehensive survey of $\gamma$-positivity in combinatorics and geometry we refer the reader to \cite{Athtba}.

Given a partition of $n$, written as $\lambda \vdash n$, we will denote by $\SYT(\lambda)$ the set of all standard Young tableaux of shape $\lambda$. The descent set $\Des(Q)$ of a standard Young tableau $Q \in \SYT(\lambda)$, where $\lambda \vdash n$, is the set of all $i \in [n-1]$ for which $i+1$ appears in a lower row in $Q$ than $i$ does. We recall that the Robinson-Schensted correspondence \cite[Theorem~7.13.5]{StaEC2} is a bijection from $\mathfrak{S}_{n}$ to the set of pairs $(P,Q)$ of standard Young tableaux of the same shape and size $n$ with the property \cite[Lemma~7.23.1]{StaEC2} that $\Des(w) = \Des(Q(w))$ and $P(w^{-1}) =Q(w),$ for all $w \in \mathfrak{S}_{n}$.  So, restricting ourselves to $I_{n}$ the Robinson-Schensted correspondence is a bijection from $I_{n}$ to the set $\SYT_{n}$ of all standard Young tableaux of size $n$ which preserves the descent set.  Strehl \cite{Str81} noticed that the map $Q \mapsto Q^{t}$, where $Q^{t}$ is the transpose tableau of $Q$, i.e. the tableau whose rows coincide with the columns of $Q$, is a bijection from the set of $Q \in \SYT_{n}$ with $\des(Q) = k$ to the set of $Q \in \SYT_{n}$ with $\des(Q) = n-1-k$, thus proving that the polynomial $I_{n}(x)$ is symmetric with center of symmetry at $(n-1)/2$.

\subsection{Symmetric and quasisymmetric functions}
\label{subsec:sym}
Our notation concerning these topics follows that of \cite[Chapter~7]{StaEC2} . Any unexplained terminology can be found there. Let $\mathbf{x} = (x_{1}, x_{2}, \dots)$ and $\mathbf{y} = (y_{1}, y_{2}, \dots)$ be sequences of pairwise commuting indeterminates. We will denote by $\Lambda^{n}$ (respectively, $\QSym^{n}$) the $\CC$-vector space of homogeneous symmetric (respectively, quasisymmetric) functions of degree $n$ in $\mathbf{x}$. The fundamental quasisymmetric function associated to $S \subseteq [n-1]$ is defined as 

\begin{equation*}
 F_{n,S}(\mathbf{x}) = \sum_{\substack{i_{1} \leq i_{2} \leq \cdots \leq i_{n} \\ j \in S \Rightarrow i_{j} < i_{j+1}}}x_{i_{1}}x_{i_{2}}\cdots x_{i_{n}}.
 \end{equation*}
 
\medskip
\noindent
The set $\{F_{n,S}(\mathbf{x}) : S \subseteq [n-1]\}$ is known to be a basis of $\QSym^{n}$.  The following well-known proposition expresses the Schur function $s_{\lambda}(\mathbf{x})$ associated to $\lambda \vdash n$ in terms of fundamental quasisymmetric functions.

\begin{proposition}{\rm(\cite[Theorem~7.19.7]{StaEC2})}\label{pro:schurtofun}
For every $\lambda \vdash n$,
\begin{equation}\label{eq:schurtofun}
s_{\lambda}(\mathbf{x}) = \sum_{Q \in \SYT(\lambda)} F_{n,\Des(Q)}(\mathbf{x}).
\end{equation}
\end{proposition}

\medskip
A consequence of the symmetry of the Robinson-Schensted correspondence is the following Cauchy-type identity for Schur functions, which will be useful in Section \ref{sec:gen}.

\begin{proposition}{\rm (\cite[Corollary~7.13.9]{StaEC2})}\label{pro:cauchytypeiden}
Let $\Par \coloneqq \cup_{n \geq 0}\Par(n)$, where $\Par(n)$ consists of all partitions of $n$. Then 
\begin{equation}\label{eq:cauchytypeiden}
\sum_{\lambda \in \Par} s_{\lambda}(\mathbf{x}) = \frac{1}{\prod_{i \geq 1}(1-x_{i})\prod_{1 \leq i < j}(1-x_{i}x_{j})}.
\end{equation}
\end{proposition}

\medskip
\subsection{Signed permutations, bitableaux and signed quasisymmetric functions}
\label{subsec:hyper}

Our notation concerning these topics mostly follows that of \cite[Section~2]{AAER17}.

The hyperoctahedral group $B_{n}$ consists of all signed permutations of length $n$, i.e. bijective maps $w : \Omega_{n} \rightarrow \Omega_{n}$ such that $w(a) = b $ implies $w(-a) = -b$ for every $a \in \Omega_{n}$. In this paper we view the hyperoctahedral group as a colored permutation group and it is convenient to use the total order 
\[ -1 <_{r} -2 <_{r} \cdots <_{r}  0 <_{r} 1 <_{r} 2 <_{r} \cdots\]
on $\ZZ$. For $w \in B_{n}$, we define 

\begin{equation*}
\des_{B}(w) = |\{ i \in \{0,1, \dots, n-1\} : w(i) >_{r} w(i+1)\}|
\end{equation*}

\noindent
and

\begin{equation*}
 \des^{B}(w) = |\{i \in \{0,1, \dots, n-1\} : w(i) > w(i+1)\}|,
\end{equation*}

\medskip
\noindent
where $w(0)\coloneqq0$ . Such an index is called a descent of $w$. In the first case, $\des_{B}$ coincides with the notion of descent for colored permutations \cite[Section~2]{Ste94} and in the second case, $\des^{B}$ coincides with the notion of descent in Coxeter groups \cite[Section~13.1]{PetEN}.
The $B_{n}$-Eulerian polynomial is defined by the formula 
\[ B_{n}(x) = \sum_{w \in B_{n}} x^{\des_{B}(w)} = \sum_{w \in B_{n}} x^{\des^{B}(x)}\]
for every positive integer $n$. $B_{n}$-Eulerian polynomials share similar properties with the classic Eulerian polynomials. In fact, they are $\gamma$-positive (thus unimodal) (see \cite[Theorem~2.10]{Athtba} for combinatorial interpretations of the corresponding $\gamma$-coefficients).

The signed descent set of $w \in B_{n}$, denoted  $\sDes(w),$ is defined as the pair $(\Des(w), \epsilon)$, where $\epsilon = (\epsilon_{1}, \epsilon_{2}, \dots, \epsilon_{n}) \in \{-,+\}^{n}$ is the sign vector with $i$th coordinate equal to the sign of $w(i)$ and $\Des(w)$ consists of the indices $i \in [n-1]$ for which either $\epsilon_{i} = +$ and $\epsilon_{i+1} = -$, or $\epsilon_{i} = \epsilon_{i+1}$ and $|w(i)| > |w(i+1)|$. For $w \in B_{n}$ we have $\des_{B}(w) = |\Des(w)|,$ if $\epsilon_{1} =+$ and $\des_{B}(w) = |\Des(w)| +1$, if $\epsilon_{1} = -$.

A bipartition of a positive integer $n$, written $(\lambda, \mu) \vdash n$, is any pair $(\lambda, \mu)$ of integer partitions of total sum $n$. A standard Young bitableau of shape $(\lambda,\mu) \vdash n$ and size $n$ is any pair $Q = (Q^{+}, Q^{-})$ of  tableaux which are strictly increasing along rows and columns, such that $Q^{+}$ has shape $\lambda$, $Q^{-}$ has shape $\mu$ and every element of $[n]$ appears exactly once as an entry of $Q^{+}$ or $Q^{-}$. The tableaux $Q^{+}$ and $Q^{-}$ are called the parts of $Q$. We will denote by $\SYT(\lambda, \mu)$ the set of all standard Young bitableaux of shape $(\lambda, \mu)$. The signed descent set  of $Q \in \SYT(\lambda,\mu)$, denoted $\sDes(Q)$, is defined as the pair $(\Des(Q), \epsilon)$, where $\epsilon = (\epsilon_{1}, \epsilon_{2},\dots, \epsilon_{n}) \in \{-,+\}^{n}$ is the sign vector with $i$th coordinate equal to the sign of the part of $Q$ in which $i$ appears and $\Des(Q)$ is the set of indices $i \in [n-1]$ for which either $\epsilon_{i} = +$ and $\epsilon_{i+1} = -$, or $\epsilon_{i} = \epsilon_{i+1} $ and $i+1$ appears in $Q$ in a lower row than $i$. For $Q \in \SYT(\lambda, \mu)$ we let $\des_{B}(Q) = |\Des(Q)|,$ if $\epsilon_{1}=+$ and $\des_{B}(Q) = |\Des(Q)| +1,$ if $\epsilon_{1} = -$. 

\begin{remark}\rm 
Our definitions of the signed descent set of a signed permutation and a standard Young bitableau are slightly different from, but equivalent to, the ones given in \cite[Definitions~2.2~and~2.3]{AAER17}.
\end{remark}

The Robinson-Schensted correspondence has a natural $B_{n}$-analogue. The Robinson-Schensted correspondence of type $B$, as described in \cite[Section~6]{Sta81} and \cite[Section~5]{AAER17}, is a bijection from $B_{n}$ to the set of pairs $(P^{B},Q^{B})$ of standard Young bitableaux of the same shape and size $n$ such that \cite[Proposition~5.1]{AAER17} $P^{B}(w^{-1}) = Q^{B}(w)$ and $\sDes(w) = \sDes(Q^{B}(w)),$ for all $w \in B_{n}$. So, restricting ourselves to signed involutions, the Robinson-Schensted correspondence of type $B$ is a bijection from $I_{n}^{B}$ to the set $\SYB_{n}$ of all standard Young bitableaux of size $n$ which preserves the signed descent set.

\begin{proposition}\label{pro:bpalindromic}
The polynomial $I_{n}^{B}(x)$ is symmetric with center of symmetry at $n/2$.
\end{proposition}

\noindent
\emph{Proof.} It suffices to describe a bijection from the set of all standard Young bitableaux $Q$ with $\des_{B}(Q) = k$ to the set of all standard Young bitableaux $Q$ with $\des_{B}(Q) = n-k$. Given a standard Young bitableau $Q = (Q^{+}, Q^{-}),$ we define its transpose $Q^{t} = ((Q^{t})^{+}, (Q^{t})^{-})$ by setting $(Q^{t})^{+} \coloneqq (Q^{-})^{t}$ and $(Q^{t})^{-} \coloneqq (Q^{+})^{t}$. We leave it to the reader to verify that the map $Q \mapsto Q^{t}$ has the required properties.
\qed

\bigskip
Different $B_{n}$-analogues of quasisymmetric functions have been suggested. The $B_{n}$-analogue of the fundamental quasisymmetric functions that we use is the one introduced by Poirier \cite[Section~3]{Poi98}. For $w \in B_{n}$ let

\begin{equation}\label{eq:bquasis}
F_{w}(\mathbf{x}, \mathbf{y}) = \sum_{\substack{i_{1}\leq i_{2}\leq \cdots \leq i_{n} \\ j \in \Des(w) \Rightarrow i_{j} < i_{j+1}}}z_{i_{1}}z_{i_{2}}\cdots z_{i_{n}}
\end{equation}

\medskip
\noindent
where $z_{i_{j}} = x_{i_{j}}$ if $\epsilon_{j} = +$ and $z_{i_{j}} = y_{i_{j}}$ if $\epsilon_{j} = -$. For a standard Young bitableau $Q$, we can similarly define $F_{Q}(\mathbf{x}, \mathbf{y})$ as in (\ref{eq:bquasis}) with $w$ replaced by $Q$. The following $B_{n}$-analogue of Proposition \ref{pro:schurtofun}, proved by Adin et al. \cite[Proposition~4.2]{AAER17}, plays a key role to the proof of Theorem \ref{thm:genfuninvoB}.

\begin{proposition}{\rm (\cite[Proposition~4.2]{AAER17})}\label{pro:bschurtofun}
For all partitions $\lambda, \mu$ 
\begin{equation}\label{eq:bschurtofun}
 s_{\lambda}(\mathbf{x})s_{\mu}(\mathbf{y}) = \sum_{Q \in \SYB(\lambda,\mu)} F_{Q}(\mathbf{x}, \mathbf{y}).
 \end{equation}
\end{proposition}

\section{Generating function for $I_{n}^{B}(x)$}
\label{sec:gen}

This section provides a proof of (\ref{eq:genfuninvo}) and proves Theorem \ref{thm:genfuninvoB}. Then, using Theorem \ref{thm:genfuninvoB}, it derives a linear recurrence formula for the coefficients of $I_{n}^{B}(x)$.

\medskip
\noindent
\emph{Proof of (\ref{eq:genfuninvo}).} For a power series $f(\mathbf{x})$, $m \in \ZZ_{>0}$ and indeterminate $t$ we write 
\[ f(t^{m}) = f(x_{1} = x_{2} = \cdots = x_{m} =t, x_{m+1} = x_{m+2} = \cdots = 0).\]
For $S \subseteq [n-1]$ we have (see the discussion before Proposition 7.19.12 in \cite{StaEC2}) 
\[ \sum_{m \geq 1} F_{n,S}(1^{m})\, x^{m-1} = \frac{x^{|S|}}{(1-x)^{n+1}}.\]
Thus, letting $S = \Des(w)$, we get 
\begin{equation}\label{eq:1.i}
\sum_{m \geq 1} F_{n,\Des(w)}(1^{m})\, x^{m-1} = \frac{x^{\des(w)}}{(1-x)^{n+1}}
\end{equation}
for every $w \in \mathfrak{S}_{n}$. Taking the sum over all $w \in I_{n}$, (\ref{eq:1.i}) becomes
\[ \frac{I_{n}(x)}{(1-x)^{n+1}}= \sum_{w \in I_{n}}\sum_{m \geq 1} F_{n, \Des(w)}(1^{m})\, x^{m-1} = \sum_{m \geq 1}\sum_{Q \in \SYT_{n}} F_{n, \Des(Q)}(1^{m})\, x^{m-1}\]
which by (\ref{eq:schurtofun}) becomes 
\[\frac{I_{n}(x)}{(1-x)^{n+1}}= \sum_{m \geq 1} \sum_{\lambda \vdash n} s_{\lambda}(1^{m})\, x^{m-1}.\]
Then, using (\ref{eq:cauchytypeiden}), we get

\begin{align*}
\sum_{n \geq 0} \, \frac{I_{n}(x)}{(1-x)^{n+1}} \, t^{n} &= \sum_{m \geq 1} \sum_{\lambda \in \Par} s_{\lambda}(1^{m})\, x^{m-1}t^{|\lambda|} \\
 &= \sum_{m\geq 1}\sum_{\lambda \in \Par} s_{\lambda}(t^{m})\, x^{m-1} \\ 
 &=\sum_{m\geq 1}\, \frac{x^{m-1}}{(1-t)^{m}(1-t^{2})^{\binom{m}{2}}}.
 \end{align*}
\qed

\bigskip
The following lemma, stated without proof as Equation (43) in \cite{Athtba}, will be used in the proof of Theorem \ref{thm:genfuninvoB}.  For a power series $f(\mathbf{x}, \mathbf{y})$, $m \in \ZZ_{>0}$ and indeterminates $p,q$ we write $f(p^{m},q^{m})$, where $x_{1} = x_{2} = \cdots = x_{m} =p,$  $x_{m+1} = x_{m+2} = \cdots = 0$ and $y_{1}=y_{2} = \cdots =y_{m} = q,$ $y_{m+1} = y_{m+2} = \cdots = 0$.

\begin{lemma} For $w \in B_{n}$ we have
\begin{equation}\label{eq:lemma.i}
\sum_{m \geq 1} F_{w}(1^{m},01^{m-1}) \, x^{m-1} = \frac{x^{\des_{B}(w)}}{(1-x)^{n+1}}.
\end{equation}
\end{lemma}

\noindent
\emph{Proof.}
Let $w \in B_{n}$. First suppose that $w(1)>0$. Then 
\[ F_{w}(1^{m},01^{m-1}) = \sum_{\substack{1\leq i_{1} \leq i_{2} \leq \cdots \leq i_{n} \leq m \\ j \in \Des(w) \Rightarrow i_{j} <i_{j+1}}}1 = \sum_{1\leq i_{1}' \leq i_{2}' \leq \cdots \leq i_{n}' \leq m - |\Des(w)|} 1.\]
But, in this case $\des_{B}(w) = |\Des(w)|$. So
 \[ F_{w}(1^{m},01^{m-1}) = \sum_{1\leq i_{1}' \leq i_{2}' \leq \cdots \leq i_{n}' \leq m - \des_{B}(w)} 1 = [x^{m-1}] \, \frac{x^{\des_{B}(w)}}{(1-x)^{n+1}},\]
where $[x^{n}]f(x)$ is the coefficient of $x^{n}$ in the formal power series $f(x),$ and (\ref{eq:lemma.i}) follows. Now suppose $w(1)<0$. Then 
\[ F_{w}(1^{m},01^{m-1}) = \sum_{\substack{1< i_{1} \leq i_{2} \leq \cdots \leq i_{n} \leq m \\ j \in \Des(w) \Rightarrow i_{j} <i_{j+1}}}1 = \sum_{1\leq i_{1}' \leq i_{2}' \leq \cdots \leq i_{n}' \leq m - |Des(w)|-1} 1.\]
But, in this case, $\des_{B}(w) = |\Des(w)|+1$ and so
 \[ F_{w}(1^{m},01^{m-1}) = \sum_{1\leq i_{1}' \leq i_{2}' \leq \cdots \leq i_{n}' \leq m - \des_{B}(w)} 1 = [x^{m-1}] \, \frac{x^{\des_{B}(w)}}{(1-x)^{n+1}}\]
 and (\ref{eq:lemma.i}) follows.
\qed

\bigskip
\noindent
\emph{Proof of Theorem \ref{thm:genfuninvoB}.}  Taking the sum over all $w \in I_{n}^{B}$, (\ref{eq:lemma.i}) becomes
\[ \frac{I_{n}^{B}(x)}{(1-x)^{n+1}} =\sum_{w \in I_{n}^{B}} \sum_{m\geq 1} F_{w}(1^{m},01^{m-1}) \, x^{m-1} = \sum_{m \geq 1} \sum_{Q \in \SYB_{n}}F_{Q}(1^{m},01^{m-1}) \, x^{m-1}.\] 
By (\ref{eq:bschurtofun}), we have 

\begin{equation*}
\frac{I_{n}^{B}(x)}{(1-x)^{n+1}} = \sum_{m \geq 1} \sum_{(\lambda, \mu) \vdash n } s_{\lambda}(1^{m})s_{\mu}(1^{m-1})\, x^{m-1}
\end{equation*}

\medskip
\noindent
and using (\ref{eq:cauchytypeiden}) twice we conclude that
\begin{align*}
\sum_{n\geq 0} \, \frac{I_{n}^{B}(x)}{(1-x)^{n+1}} \, t^{n} &= \sum_{m \geq 1} \sum_{\lambda, \mu \in \Par} s_{\lambda}(1^{m})s_{\mu}(1^{m-1})\, x^{m-1}\, t^{|\lambda|+|\mu|} \\
&= \sum_{m \geq 1} \left( \sum_{\lambda \in \Par}s_{\lambda}(t^{m})\right)\left(\sum_{\mu \in \Par} s_{\mu}(t^{m-1})\right) \\
&=\sum_{m \geq 1} \frac{x^{m-1}}{(1-t)^{m}(1-t^{2})^{\binom{m}{2}}(1-t)^{m-1}(1-t^{2})^{\binom{m-1}{2}}} \\
&= \sum_{m \geq 1} \frac{x^{m-1}}{(1-t)^{2m-1}(1-t^{2})^{(m-1)^{2}}}.
\end{align*}
\qed

\bigskip
Let $I_{n}^{B}(x) = \sum_{k=0}^{n}I_{n,k}^{B}x^{k}$. The following proposition provides a linear recurrence formula for the sequence $I_{n,0}^{B}, I_{n,1}^{B}, \dots, I_{n,n}^{B}$ for every $n \in \ZZ_{>0}$. For a similar formula regarding the coefficients of $I_{n}(x),$ see \cite[Theorem~2.2]{GZ06}.

\begin{proposition}\label{pro:formula}
For $n \geq 3$ and $k \geq 0$, the numbers $I_{n,k}^{B}$ satisfy the following recurrence formula: 
 \begin{align}\label{eq:formula}
 nI_{n,k}^{B} =& \; (2k+1)I_{n-1,k}^{B} + (2n-2k+1)I_{n-1,k-1}^{B} + (n-1 + 2k(k+1))I_{n-2,k}^{B} \nonumber \\
 &+ (2(n-1) + 4(n-k-1)(k-1))I_{n-2, k-1}^{B} \nonumber \\
  &+ ((2n-3)(n-1) + 2(k-2)(k-2n+1))I_{n-2,k-2}^{B} 
 \end{align}
 where $I_{n,k}^{B} =0$ for every $k <0$.
\end{proposition}

\noindent
\emph{Proof.} Extracting the coefficients of $t^{n}$ of both sides in (\ref{eq:genfuninvoB}) we get
\begin{equation}\label{eq:pro.i}
 \frac{I_{n}^{B}(x)}{(1-x)^{n+1}} = \sum_{m \geq 0}\sum_{j = 0}^{\lfloor n/2 \rfloor} {{m^{2}+j-1}\choose{j}}{{2m+n -2j}\choose{n-2j}}x^{m}. 
\end{equation}
Let 
\[ r(n,m) = \sum_{j = 0}^{\lfloor n/2 \rfloor} {{m^{2}+j-1}\choose{j}}{{2m+n -2j}\choose{n-2j}}\]
for all $n,m \geq 0$. Clearly, $\sum_{n \geq 0}r(n,m)t^{n} = (1-t)^{-(2m+1)}(1-t^{2})^{-m^{2}}$ and taking derivatives with respect to $t$ we have

\begin{equation*}
 \sum_{n \geq 1}nr(n,m) t^{n-1} = \frac{(2m^{2} + 2m +1)t +2m +1}{1-t^{2}}\sum_{n \geq 0}r(n,m) t^{n}.
 \end{equation*}
 
\medskip
\noindent
Extracting the coefficients of $t^{n}$ in the expression above we obtain the following recurrence formula for the numbers $r(n,m)$:
\begin{equation}\label{eq:helpformula}
nr(n,m) = (2m+1)r(n-1,m) + (2m^{2} + 2m + n -1)r(n-2,m),
\end{equation}
for every $n \geq 3$. Substituting (\ref{eq:helpformula}) to (\ref{eq:pro.i})
we have
\begin{align*} \frac{nI_{n}^{B}(x)}{(1-x)^{n+1}} =& \; 2x\left(\frac{I_{n}^{B}(x)}{(1-x)^{n+1}}\right)' + \frac{I_{n-1}^{B}(x)}{(1-x)^{n}} + 2x^{2} \left(\frac{I_{n-2}^{B}(x)}{(1-x)^{n-1}}\right)'' \\
&+ 4x\left(\frac{I_{n-2}^{B}(x)}{(1-x)^{n-1}}\right)' + (n-2) \frac{I_{n-2}^{B}(x)}{(1-x)^{n-1}}
\end{align*}
or
\begin{align*}
nI_{n}^{B}(x) =& \; (2n-1)xI_{n-1}^{B}(x) + I_{n-1}^{B}(x) -2x^{2}I_{n-1}^{B'}(x) + 2xI_{n-1}^{B'}(x) \\
&+ (2n-3)(n-1)x^{2}I_{n-2}^{B}(x) + 2(n-1)xI_{n-2}^{B}(x) + (n-1)I_{n-2}^{B}(x) \\
&-  4(n-2)x^{3}I_{n-2}^{B'}(x) + 4(n-3)x^{2}I_{n-2}^{B'}(x) + 4xI_{n-2}^{B'}(x) \\
&+ 2x^{4}I_{n-2}^{B''}(x) - 4x^{3}I_{n-2}^{B''}(x) + 2x^{2}I_{n-2}^{B''}(x).
\end{align*}
Comparing the coefficients of $x^{k}$ in both sides of the above identity, we obtain the desired formula (\ref{eq:formula}).
\qed

\begin{remark} \rm 
Another proof of Proposition \ref{pro:bpalindromic} can be obtained from Proposition \ref{pro:formula} as follows. The right-hand side of (\ref{eq:formula}) is invariant under the substitution $k \rightarrow n-k$ provided the sequences $I_{n-1,k}^{B}$ and $I_{n-2,k}^{B}$ are symmetric. Thus, by induction the proof follows.
\end{remark}
\section{Unimodality of $I_{n}^{B}(x)$}
\label{sec:uni}

This section proves the main result of this paper, namely the unimodality of $I_{n}^{B}(x)$  for every positive integer $n$. The method is similar to that of Guo and Zeng \cite{GZ06} who proved  the unimodality of $I_{n}(x)$. First, we recall the following observation of Guo and Zeng \cite[Lemma~3.1]{GZ06}.

\begin{lemma}\label{lemma}
Let $x_{0}, x_{1}, \dots, x_{n}$ and $a_{0}, a_{1}, \dots, a_{n}$ be real numbers such that $x_{0} \geq x_{1} \geq \cdots \geq x_{n} \geq 0$ and $a_{0} + a_{1} + \cdots + a_{n} \geq 0$ for all $0 \leq k \leq n$. Then 
\[ a_{0}x_{1} + a_{1}x_{1} + \cdots + a_{n}x_{n} \geq 0.\]
\end{lemma}

By Proposition \ref{pro:bpalindromic} the sequence $I_{n,0}^{B}, I_{n,1}^{B}, \dots ,I_{n,n}^{B}$ is symmetric. So, to prove its unimodality, it suffices to show that $I_{n,k}^{B} \geq I_{n,k-1}^{B}$ for every $0 \leq k \leq \lfloor n/2 \rfloor$.

\begin{theorem}\label{thm:unimodal}
The sequence $I_{n,0}^{B}, I_{n,1}^{B}, \dots, I_{n,n}^{B}$ is unimodal for every $n \in \ZZ_{>0}$.
\end{theorem}

\noindent
\emph{Proof.} We proceed by induction on $n$. For $n \leq 2$ we saw in Section \ref{sec:intro} that the statement is true. Suppose that $n>2$ and that the sequences $I_{n-1,k}^{B}$ and  $I_{n-2,k}^{B}$ are unimodal. Substituting $k \rightarrow k-1$ in (\ref{eq:formula}) we obtain

\begin{align}\label{eq:thm.i}
nI_{n,k-1}^{B} =& \; (2k-1)I_{n-1,k-1}^{B} + (2n -2k +3)I_{n-1, k-2}^{B} + (n-1 + 2k(k-1))I_{n-2,k-1}^{B} \nonumber \\
&+ (2(n-1) + 4(n-k)(k-2))I_{n-2,k-2}^{B} \nonumber \\ 
&+ ((2n-3)(n-1) + 2(k-3)(k-2n))I_{n-2,k-3}^{B}.
\end{align}
Subtracting (\ref{eq:thm.i}) from (\ref{eq:formula}) we get
\begin{align}
n(I_{n,k}^{B} - I_{n,k-1}^{B}) =& \; A_{0}I_{n-1,k}^{B} + A_{1}I_{n-1,k-1}^{B} + A_{2}I_{n-1,k-2}^{B} \nonumber \\
&+ D_{0}I_{n-2,k}^{B} + D_{1}I_{n-2,k-1}^{B} + D_{2}I_{n-2,k-2}^{B} + D_{3}I_{n-2,k-3}^{B},
\end{align}
where
\begin{align*}
A_{0} &= 2k +1 & D_{1} &= 4nk -3n -6k^{2} + 2k + 3 \\
A_{1} &= 2n-4k+2  & D_{2} &= 2n^{2} - 8nk + 9n + 6k^{2} - 10k +1 \\
A_{2} &= -2n + 2k -3  & D_{3} &= -2n^{2} + 4nk -7n -2k^{2} + 6k -3.\\
D_{0} &= n + 2k^{2} + 2k -1 &
\end{align*}

By the induction hypothesis we know that $I_{n-1,k}^{B} \geq I_{n-1,k-1}^{B} \geq I_{n-1,k-2}^{B}$ for every $0 \leq k \leq \lfloor n/2 \rfloor$. Also, $A_{0} = 2k +1 \geq 0$, $A_{0} + A_{1} = 2(n-k) +3 \geq 0$ and $A_{0} + A_{1} + A_{2} =0$. Thus, by Lemma \ref{lemma}, it follows that
\begin{equation}\label{eq:part1}
 A_{0}I_{n-1,k}^{B} + A_{1}I_{n-1,k-1}^{B} + A_{2}I_{n-1,k-2}^{B} \geq 0
\end{equation}
for every $0 \leq k \leq \lfloor n/2 \rfloor$. It remains to show that 
\begin{equation}\label{eq:part2}
D_{0}I_{n-2,k}^{B} + D_{1}I_{n-2,k-1}^{B} + D_{2}I_{n-2,k-2}^{B} + D_{3}I_{n-2,k-3}^{B} \geq 0.
\end{equation}
Indeed, combining (\ref{eq:part1}) and (\ref{eq:part2}) yields $I_{n,k}^{B} - I_{n,k-1}^{B} \geq 0$ for every $0 \leq k \leq \lfloor n/2 \rfloor$ and the proof is completed.

We distinguish two cases. In the first case, we assume that $1 \leq k \leq \lfloor n/2 \rfloor -1$. By the induction hypothesis we know that $I_{n-2,k}^{B} \geq I_{n-2,k-1}^{B} \geq I_{n-2,k-2}^{B} \geq I_{n-2,k-3}^{B}$. Also, 
\begin{align*}
D_{0} &= n + 2k^{2} + 2k -1 \geq 0 \\
D_{0} + D_{1} &= 2(k-1)(n-2k) + 2 \geq 0 \\
D_{0} + D_{1} + D_{2} &= 2n(n-2k) + n + 6(n-k) + 2k^{2} + 3 \geq 0 \\
D_{0} + D_{1} + D_{2} + D_{3} &= 0.
\end{align*}
Thus, by Lemma \ref{lemma}, (\ref{eq:part2}) follows. For the second case, suppose that $k = \lfloor n/2 \rfloor$. If $n$ is odd, then by the induction hypothesis and symmetry we have

\begin{equation*}
I_{n-2,k}^{B} = I_{n-2,k-1}^{B} \geq I_{n-2,k-2}^{B} \geq I_{n-2,k-3}^{B}
\end{equation*}

\medskip
\noindent
and from the previous calculations (\ref{eq:part2}) follows.  If $n$ is even, then by the induction hypothesis and symmetry we have

\begin{equation*}
I_{n-2,k}^{B} \geq  I_{n-2,k-1}^{B} = I_{n-2,k-2}^{B} \geq I_{n-2,k-3}^{B}.
\end{equation*}

\medskip
\noindent
In this case, we notice that $D_{1} = \frac{1}{2}(n^{2} -4n +6) = \frac{1}{2}((n-2)^{2} + 2)>0$. Therefore, by Lemma \ref{lemma} and the previous calculations, (\ref{eq:part2}) follows. \qed

\section{Remarks and open problems}
\label{sec:remarks}

This section discusses further properties of $I_{n}^{B}(x)$. Barnabei et al. \cite[Theorem~7]{BBS09} proved that $I_{n}(x)$ is not log-concave, answering negatively a conjecture due to Brenti. The same holds for $I_{n}^{B}(x)$.
 
\begin{proposition}\label{pro:blogconcavity}
The sequence $I_{n,0}^{B}, I_{n,1}^{B}, \dots, I_{n,n}^{B}$ is not log-concave in general.
\end{proposition}

\medskip
\noindent
\emph{Proof.} 
It is known that the product of a unimodal and log-concave polynomial with a log-concave polynomial is log-concave as well (see \cite[Proposition~6]{BBS09}). Formula \ref{eq:pro.i} shows that the polynomial $r(x) = \sum_{k=0}^{n}r(n,k)x^{k}$ is equal to the product of $I_{n}^{B}(x)$ and $q(x) = \sum_{k=0}^{n} {{n+k}\choose{k}}x^{k}$. Clearly, $q(x)$ is log-concave. So, if $I_{n}^{B}(x)$ was log-concave then $r(x)$ has to be log-concave since $I_{n}^{B}(x)$ is unimodal by Theorem \ref{thm:unimodal}. But, $r(x)$ is not log-concave in general. For instance, 
\[ r(89,2)^{2} = 113 789 153 706 560 010 000 < 114 890 217 312 335 629 500 = r(89,1)r(89,3).\]\qed 

\medskip
Guo and Zeng \cite[Conjecture~4.1]{GZ06} conjectured that $I_{n}(x)$ is $\gamma$-positive. The conjecture remains unsolved as far as we know. We notice that 
\[ I_{n}^{B}(x) =\begin{cases}
 1+x, & \text{if } n=1  \\
 (1+x)^{2} + 2x, & \text{if } n=2 \\
 (1+x)^{3} + 6x(1+x), & \text{if } n=3 \\
 (1+x)^{4} + 13x(1+x)^{2} + 8x^{2}, & \text{if } n=4\\
 (1+x)^{5} + 23x(1+x)^{3} + 48x^{2}(1+x), & \text{if } n=5\\
(1+x)^{6} + 37x(1+x)^{4} + 168x^{2}(1+x)^{2} + 56x^{3}, & \text{if } n=6.
\end{cases}\]
The above calculations and the similarity between $I_{n}(x)$ and $I_{n}^{B}(x)$ suggest that $I_{n}^{B}(x)$ may be $\gamma$-positive for every $n$ (see also \cite[Section~2.1.5]{Athtba}).

\begin{remark} \rm
As mentioned in \cite[Section~2.1.5]{Athtba}, the definition of $I_{n}^{B}(x)$ is presumably but not obviously unaffected when $\des_{B}$ is replaced by $\des^{B}$. We have confirmed it for $n \leq 5$.
\end{remark}

\section*{Acknowledgments}
The author would like to thank Christos Athanasiadis for suggesting the problem of unimodality of $I_{n}^{B}(x)$ and for valuable discussions. This work is part of the author's Master thesis at the Department of Mathematics of the National and Kapodistrian University of Athens under the supervision of Athanasiadis.

\end{document}